\documentclass[letterpaper,10pt,conference]{ieeeconf}%

\usepackage{hyperref} 
\usepackage{graphicx}
\usepackage{amsthm}
\usepackage{algorithm}
\usepackage[noend]{algpseudocode}
\usepackage[cmex10]{amsmath}
\usepackage{amsfonts,amssymb}
\usepackage{times}
\usepackage{bm}
\usepackage{bbm}
\usepackage{color}
\usepackage{graphicx}
\usepackage{subcaption}
\usepackage{amsfonts}
\usepackage{cite}
\usepackage{mathtools}
\usepackage{amssymb}%
\providecommand{\U}[1]{\protect\rule{.1in}{.1in}}

\IEEEoverridecommandlockouts
\newtheorem{theorem}{Theorem}

\newtheorem{lemma}{Lemma}
\newtheorem{corollary}{Corollary}
\newtheorem{definition}{Definition}
\theoremstyle{definition}

\makeatletter
\def\BState{\State\hskip-\ALG@thistlm}
\makeatother

\begin{document}

\title{{\LARGE \textbf{Decentralized Event-Driven Algorithms for Multi-Agent
Persistent Monitoring }}}
\author{Nan Zhou$^{1}$, Christos G. Cassandras$^{1,2}$, Xi Yu$^{3}$, and Sean B. Andersson$^{1,3}$ \\{\small $^{1}$Division of Systems Engineering, $^{2}$Department of Electrical and Computer Engineering, $^{3}$Department of Mechanical
Engineering}\\Boston University, Boston, MA 02215, USA\\E-mail:\texttt{\{nanzhou,cgc,xyu,sanderss\}@bu.edu \thanks{* The work of Cassandras and Zhou is supported in part by NSF under grants ECCS-1509084, IIP-1430145, and CNS-1645681, by AFOSR under grant FA9550-15-1-0471, by DOE under grant 46100 and by the MathWorks. The work of Andersson and Yu is supported in part by the NSF through grant
ECCS-1509084 and CMMI-1562031.} }}
\maketitle

\begin{abstract}
We address the issue of identifying conditions under which the centralized solution to the optimal multi-agent persistent monitoring problem can be recovered in a decentralized event-driven manner. In this problem, multiple agents interact with a finite number of targets and the objective is to control their movement in order to minimize an uncertainty metric associated with the targets. In a one-dimensional setting, it has been shown that the optimal solution can be reduced to a simpler parametric optimization problem and that the behavior of agents under optimal control is described by a hybrid system. This hybrid system can be analyzed using Infinitesimal Perturbation Analysis (IPA) to obtain a complete on-line solution through an event-driven centralized gradient-based algorithm. We show that the IPA gradient can be recovered in a distributed manner in which each agent optimizes its trajectory based on local information, except for one event requiring communication from a non-neighbor agent. Simulation examples are included to illustrate the effectiveness of this \textquotedblleft almost decentralized\textquotedblright\ algorithm and its fully decentralized counterpart where the aforementioned non-local event is ignored.

\end{abstract}


\section{Introduction}

\label{sec:intro} Systems consisting of cooperating mobile agents are often
used to perform tasks such as coverage \cite{zhong2011distributed},
surveillance \cite{michael2011persistent}, or environmental sampling
\cite{leonard2010coordinated}. A \emph{persistent monitoring} task is one
where agents must cooperatively monitor a dynamically changing environment
that cannot be fully covered by a stationary team of agents (as in coverage
control) \cite{cassandras2013optimal}. Once the exploration process leads to
the discovery of various \textquotedblleft points of
interest\textquotedblright, then these become \textquotedblleft data
sources\textquotedblright\ or \textquotedblleft targets\textquotedblright%
\ which need to be perpetually monitored. Thus, in contrast to sweep coverage
and patrolling where \emph{every} point in a mission space is of interest
\cite{Smith:2012fq,lin2015optimal}, the problem we address here focuses on a
\emph{finite number} of data sources or \textquotedblleft
targets\textquotedblright\ (typically larger than the number of agents).

In this setting, the agents interact with targets through their sensing
capabilities which are normally dependent upon their physical distance from
the target. The uncertainty state of a target increases when no agent is
visiting it and decreases when it is being monitored by one or more agents
(i.e., it is within their sensing range). The objective is to minimize an
overall measure of target uncertainty states by controlling the movement of
all agents in a cooperative manner. Unlike many other multi-agent systems
modeled solely through a network of interconnected agents, here we have two
networks, one whose nodes are agents and one whose nodes are targets. Since
agents interact with targets, this interaction is modeled by establishing
links between nodes belonging to the two different networks. Moreover, since
agents are mobile, the overall graph topology in such systems is time-varying.
Thus, the resulting complexity of this class of problems is significant. This
has motivated approaches where rather than viewing these as agent-to-target
assignment problems \cite{stump2011multi,yu2017optimal} (which are
computationally intensive and do not scale well in the number of targets and
agents), one treats them as trajectory design and optimization problems
\cite{cassandras2013optimal},\cite{horling2004survey}.

In \cite{zhou2016optimal}, we studied the persistent monitoring problem in a
one-dimensional (1D) mission space and showed that it can be formulated as an
optimal control problem whose solution is parametric, i.e., the optimal
control problem is reduced to a parametric optimization one. In particular,
every optimal agent trajectory is characterized by a finite number of points
where the agent switches direction and by a dwell time at each such point.
As a result, the behavior of agents under optimal control is described by a
hybrid system. This allows us to make use of Infinitesimal Perturbation
Analysis (IPA) \cite{cassandras2010perturbation,wardi2010unified} to determine
on-line the gradient of the objective function with respect to these
parameters and to obtain a (possibly local) optimal trajectory. Our approach
exploits IPA's event-driven nature to render it \emph{scalable} in the number
of \emph{events} in the system and not its state space.

The optimal controller developed in \cite{zhou2016optimal} is established
based on the assumption that agents are all connected under a centralized
controller which can provide information and coordinate all agents. Similar
centralized controllers for such problems can be found in
\cite{leahy2016provably,leonard2010coordinated,cassandras2013optimal}.
Clearly, a centralized controller can be energy-consuming due to communication
costs \cite{zhong2010asynchronous} and unreliable in adversarial environments.
In this paper, we address the question of whether it is possible to develop
\emph{decentralized} controllers for persistent monitoring problems with a
finite numbers of targets to be monitored.

Decentralization aims to achieve the same global objective as a central
controller by distributing functionality to the agents so that each one acts
based on local information or by communicating with only a set of neighbors.
Such distributed algorithms have been derived and applied to coverage control
\cite{zhong2011distributed}, formation control \cite{ren2008distributed}, and
consensus problems \cite{olfati2007consensus} where we usually assume a static
fully connected network environment. On the other hand, decentralization in a
persistent monitoring setting is particularly challenging due to the
time-varying nature of the agent network and the fact that agents take actions
depending on interactions with the environment (targets) which cannot be
easily shared through the agent network.

The contribution of this paper consists of identifying explicit conditions
under which the centralized solution to the optimal persistent monitoring
problem studied in \cite{zhou2016optimal} can be recovered through an
\textquotedblleft almost decentralized\textquotedblright\ and entirely
event-driven manner. In particular, each agent uses $(i)$ its own local
information (to be precisely defined later), $(ii)$ information (in the form
of observable events) from agents that happen to be its neighbors at the time
such events occur, and $(iii)$ a single specific event type communicated by a
non-neighbor agent when it occurs. It is the latter that prevents a completely
decentralized control scheme, although, as we will see, ignoring this
non-local event results in little loss of accuracy. In addition, we develop
such an \textquotedblleft almost decentralized\textquotedblright\ algorithm
which, compared to the centralized solution in \cite{zhou2016optimal},
significantly reduces communication costs while yielding the same performance.
The main decentralization result exploits the structure of the IPA gradient of
the objective function: the gradient component associated with an agent turns
out to depend only on a limited number of events, all of which are local or
observed by (time-varying) neighbors except for one event requiring
communication with a non-neighbor when it occurs. Moreover, this IPA gradient
structure is not limited to the 1D problem considered in this paper, but
extends to its 2D version as well, a direction we are pursuing in ongoing research.

The paper is organized as follows. In Section II we present the persistent
monitoring problem formulation and introduce different neighborhood concepts
that capture the interaction between agents and targets. In Section III, we
review the optimal control solution of the problem in the 1D setting and in
Section IV carry out IPA for the resulting hybrid system. In Section V, we
show that the optimal control solution can be decentralized with only one non-neighbor event needed by an agent to derive it
in a distributed manner. We provide simulation examples in Section VI to
illustrate the resulting algorithm and its performance.

\section{Problem formulation}

\label{sec:PM_form} We begin by reviewing the persistent monitoring model and
problem formulation introduced in \cite{zhou2016optimal}.

\textbf{Agent dynamics}. We consider $N$ agents moving in a one-dimensional
(1D) mission space $[0,L]\subset\mathbb{R}$. Each agent can control its speed
and direction. The speed input is scaled and bounded in $[-1,1]$. The position
of each agent $j$ is represented as $s_{j}(t)\in\lbrack0,L]$ with the
dynamics:
\begin{equation}
\dot{s}_{j}(t)=u_{j}(t),\quad|u_{j}(t)|\leq1,\quad\forall j=1,2\ldots,N
\label{eq:DynOfS}%
\end{equation}

\textbf{Agent sensing model}. The ability of an agent to sense its environment
is modeled by a function $p_{j}(x,s_{j})$ that measures the probability that
an event at location $x\in\left[  0,L\right]  $ is detected by agent $j$ at
$s_{j}(t)$. We assume that $p_{j}(x,s_{j})=1$ if $x=s_{j}$, and that
$p_{j}(x,s_{j})$ is monotonically non-increasing in the distance $\Vert
x-s_{j}\Vert$, thus capturing the reduced effectiveness of a sensor over its
range. We consider this range to be finite and denoted by $r_{j}$. Although
our analysis is not affected by the precise sensing model $p_{j}(x,s_{j})$, we
will limit ourselves to a linear decay model as follows:%
\begin{equation}
p_{j}(x,s_{j})=\max\left\{  0,1-\frac{\Vert x-s_{j}\Vert}{r_{j}}\right\}
\label{eq:sensing_single}%
\end{equation}
Unlike the sweep coverage problem, here we consider a known finite set of
targets located at $x_{i}\in\lbrack0,L],$ $i=1,\ldots,M$. We then set
$p_{j}(x_{i},s_{j}(t))\equiv p_{ij}(s_{j}(t))$ for simplicity. For $N$ agents
sensing simultaneously, assuming detection independence, the sensing
capability of $N$ agents on target $i$ can be captured by the joint detection
probability function
\begin{equation}
P_{i}\left(  \mathbf{s}(t)\right)  =1-\prod_{j=1}^{N}\left(  1-p_{ij}%
(s_{j}(t))\right)  \label{eq:joint_sensing}%
\end{equation}
where we set $\mathbf{s}(t)=[s_{1}\left(  t\right)  ,\ldots,s_{N}\left(
t\right)  ]^{\text{T}}$.

\textbf{Target dynamics}. We define uncertainty functions $R_{i}(t)$
associated with targets $i=1,\ldots,M$, so that they have the following
properties: $(i)$ $R_{i}(t)$ increases with a prespecified rate $A_{i}$ if
$P_{i}\left(  \mathbf{s}(t)\right)  =0$ (as shown in \cite{zhou2016optimal},
this can be allowed to be a random process $\{A_{i}(t)\}$), $(ii)$ $R_{i}(t)$
decreases with a fixed rate $B_{i}$ if $P_{i}\left(  \mathbf{s}(t)\right)  =1$
and $(iii)$ $R_{i}(t)\geq0$ for all $t$. It is then natural to model
uncertainty dynamics associated with each target as follows:
\begin{equation}
\dot{R}_{i}(t)=\hspace{-0.1cm}\left\{
\begin{array}
[c]{ll}%
0 \quad\quad\quad\text{if }R_{i}(t)=0\text{ and }A_{i}\leq
B_{i}P_{i}\left(  \mathbf{s}(t)\right)  & \\
A_{i}-B_{i}P_{i}\left(  \mathbf{s}(t)\right)  \quad\text{otherwise} &
\end{array}
\right.  \label{eq:multiDynR}%
\end{equation}
where we assume that initial conditions $R_{i}(0)$, $i=1,\ldots,M$ are given
and that $B_{i}>A_{i}>0$ to ensure a strict decrease in $R_{i}(t)$ when
$P_{i}(\mathbf{s}(t))=1$.

\textbf{Optimal control problem}. Our goal is to control the movement of the
$N$ agents through $u_{j}\left(  t\right)  $ in (\ref{eq:DynOfS}) so that the
cumulative average uncertainty over all targets $i=1,\ldots,M$ is minimized
over a fixed time horizon $T$. Thus, setting $\mathbf{u}\left(  t\right)
=\left[  u_{1}\left(  t\right)  ,\ldots,u_{N}\left(  t\right)  \right]
^{\text{T}}$ we aim to solve the following optimal control problem:
\begin{equation}
\mathbf{P1:}\quad\min_{\mathbf{u}\left(  t\right)  }\text{ \ }J=\frac{1}%
{T}\int_{0}^{T}\sum_{i=1}^{M}R_{i}(t)dt \label{eq:costfunction}%
\end{equation}
subject to the agent dynamics (\ref{eq:DynOfS}) and target uncertainty dynamics
(\ref{eq:multiDynR}). Generally, the classical solution of
(\ref{eq:costfunction}) involves solving a Two Point Boundary Value Problem
(TPBVP) which requires global information of all agents and targets. In this
paper, we will limit the information of each agent to itself and its neighbors
and study whether this objective function can be optimized in a distributed manner.

\begin{figure}[pb]
\centering
\includegraphics[width=0.9\linewidth]{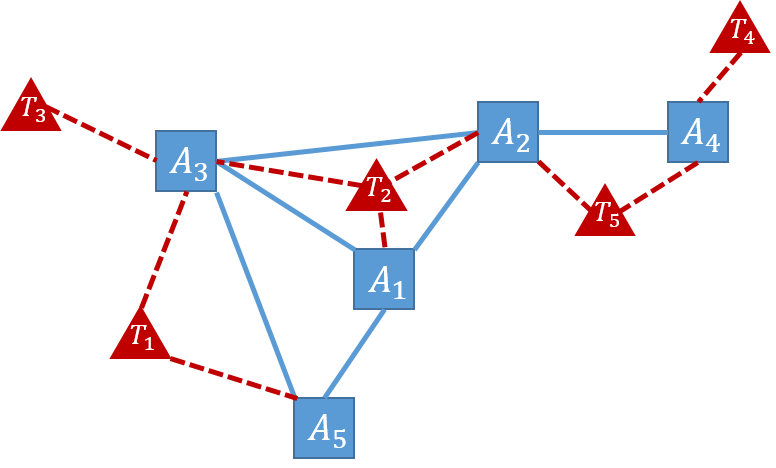}\caption{{\protect\small Agent-target
network. Red triangles are targets and blue squares are agents. Blue lines
indicate the neighbor agents of an agent and red lines indicate the neighbor
targets of an agent.}}%
\label{fig:agent_target_networks}%
\end{figure}

\textbf{Limited information model for decentralization.}
In our model, an agent is capable of observing information within its sensing
range, specifically the state $R_{i}(t)$ of all targets $i$ such that
$p_{ij}(s_{j}(t))>0$. Moreover, agents can communicate with their neighboring
agents to acquire information such as agent positions, speeds, and the states
of targets which are within their own sensing ranges. In contrast to
traditional multi-agent systems modeled through a network of agents, in the
persistent monitoring setting agents move to interact with targets as shown in
Fig. \ref{fig:agent_target_networks}. Therefore, the network model includes
both agents and targets and we need to revisit the concept of neighborhood,
accounting as well for the fact that neighborhoods are time-varying. We begin
with the observation that agents have two types of neighbors: nearby agents
and nearby targets. On the other hand, the neighborhood of a target consists
of just nearby agents. We do not explicitly model any possible connectivity
among targets; however, if the target topology is fully connected, then it is
possible for an agent near one target to acquire information about all targets.

\begin{definition}
The agent neighborhood of agent $j$ is the set $\mathcal{A}_{j}(t)=\{k:\Vert
s_{k}(t)-s_{j}(t)\Vert\leq r_{c},$ $k\neq j,$ $k=1,\ldots,N\}$.
\end{definition}

This is a conventional definition of neighbors in multi-agent systems, where
$r_{c}$ is a communication range, but we point out that it is
\emph{time-dependent} since agents are generally moving. As an example, in
Fig. \ref{fig:agent_target_networks}, $\mathcal{A}_{1}=\{A_{2},A_{3},A_{5}\}$.

\begin{definition}
\label{def:neighboring_targets} The target neighborhood of agent $j$ is the
set $\mathcal{T}_{j}(t)=\{i:|x_{i}-s_{j}(t)|\leq r_{j},$ $i=1,\ldots,M\}.$
\end{definition}

This includes all targets which are within agent $j$'s sensing range. In Fig.
\ref{fig:agent_target_networks}, $\mathcal{T}_{3}=\{T_{1},T_{2},T_{3}\}$.
Assuming the agents are homogeneous with a common sensing range $r$, we
require that $r_{c}\geq2r$ in order to establish communication among agents
that are sensing the same target.

\begin{definition}
\label{def:neighboring_agents} The agent neighborhood of target $i$ is the set
$\mathcal{B}_{i}(t)=\{j:|s_{j}(t)-x_{i}|\leq r_{j},$ $j=1,\ldots,N\}.$
\end{definition}

This set captures all the neighbor agents of target $i$. In Fig.
\ref{fig:agent_target_networks}, $\mathcal{B}_{2}=\{A_{1},A_{2},A_{3}\}$.
Using Definition \ref{def:neighboring_agents}, the joint sensing probability
in \eqref{eq:joint_sensing} can be rewritten as:
\begin{equation}
P_{i}\left(  \mathbf{s}(t)\right)  =1-\prod_{j\in\mathcal{B}_{i}(t)}\left(
1-p_{ij}(s_{j}(t))\right)  \label{JointSensingProb}%
\end{equation}
where $\mathcal{B}_{i}(t)\subseteq\{1,\ldots,N\}$. We further define
\begin{equation}
\mathcal{N}_{ij}(t)=\mathcal{B}_{i}(t)\setminus\{j\}
\label{eq:new_def_neighbor}%
\end{equation}
to indicate the \textquotedblleft collaborators\textquotedblright\ of agent
$j$ in sensing target $i$. Note that $\mathcal{N}_{ij}(t)=\{k:k\in
A_{j}(t)\text{ and }k\in B_{i}(t)\}$, thus capturing a neighbor of agent $j$
and target $i$ at the same time.

Our limited information model restricts observations of each agent to the
agent's sensing range. However, any agent $j$ is allowed to communicate with
its neighbors in $\mathcal{A}_{j}(t)$. Therefore, the local information of an
agent is the union of the observations of agent $j$ and the observations of
agents $k\in\mathcal{A}_{j}(t)$. In Section \ref{sec:Events}, we will
explicitly define the precise meaning of \textquotedblleft
information\textquotedblright\ above to consist of observable events such as
\textquotedblleft agent stops\textquotedblright\ or \textquotedblleft target
state becomes $R_{i}(t)=0$\textquotedblright. In Section \ref{sec:DecentOpt}, we will show how
\textbf{P1} can be solved by each agent under this limited information model
as opposed to the centralized one in \cite{zhou2016optimal}.

\section{From optimal control to parametric optimization}

\label{sec:Opt_control} In this section, we review properties of the
centralized optimal control solution of \textbf{P1} which allow it to be
reduced to a parametric optimization problem \cite{zhou2016optimal}. This
leads to the use of the Infinitesimal Perturbation Analysis (IPA) gradient
estimation approach \cite{cassandras2010perturbation} to find an explicit
solution through a gradient-based algorithm. We begin by defining the state
vector $\mathbf{x}(t)=[R_{1}(t),...R_{M}(t),s_{1}(t)...s_{N}(t)]$ and
associated costate vector $\lambda=[\lambda_{1}(t),...,\lambda_{M}%
(t),\lambda_{s_{1}}(t),...,\lambda_{s_{N}}(t)]$. Due to the discontinuity in
the dynamics of $R_{i}(t)$ in (\ref{eq:multiDynR}), the optimal state
trajectory may contain a boundary arc when $R_{i}(t)=0$ for some $i$;
otherwise, the state evolves in an interior arc. Using (\ref{eq:DynOfS}) and
(\ref{eq:multiDynR}), the Hamiltonian is
\begin{equation}
H(\mathbf{x},\lambda,\mathbf{u})=\sum_{i=1}^{M}R_{i}(t)+\sum_{i=1}^{M}%
\lambda_{i}(t)\dot{R}_{i}(t)+\sum_{j=1}^{N}\lambda_{s_{j}}(t)u_{j}(t)
\label{eq:Hamiltonian}%
\end{equation}
Applying the Pontryagin Minimum Principle to (\ref{eq:Hamiltonian}) with
$\mathbf{u}^{\star}(t)$, $t\in\lbrack0,T)$, denoting an optimal control, a
necessary condition for optimality is
\begin{equation}
u_{j}^{\ast}(t)=%
\begin{cases}
1 & \quad\text{if }\lambda_{s_{j}}(t)<0\\
-1 & \quad\text{if }\lambda_{s_{j}}(t)>0
\end{cases}
\label{eq:NecessaryConditionPMP}%
\end{equation}
Note that there exists a possibility that $\lambda_{s_{j}}\left(  t\right)
=0$ over some finite singular intervals \cite{bryson1975applied}, in which
case $u_{j}^{\ast}(t)$ may take values in $\{-1,0,1\}$.

\begin{figure}[h]
\centering
\includegraphics[width=0.9\linewidth]{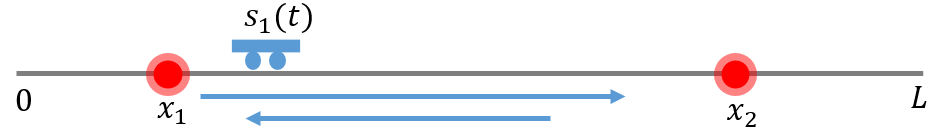}
\includegraphics[width=0.95\linewidth]{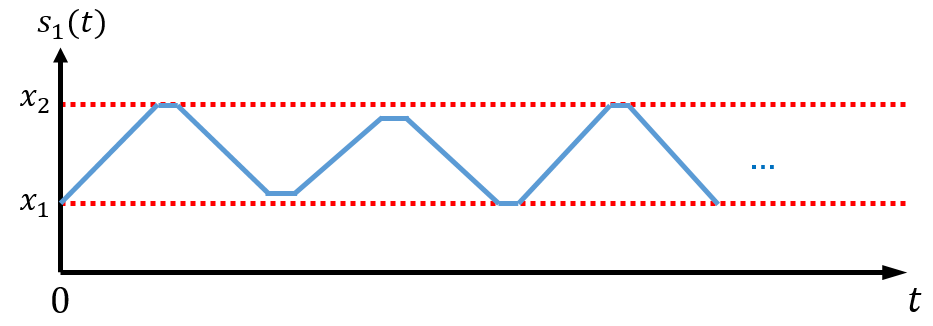}\caption{{\protect\small An
example of an agent monitoring two targets. The bottom shows the agent
trajectory which can be parameterized to a sequence of control switching
points and associated dwell times.}}%
\label{fig:para_illus}%
\end{figure}

\begin{figure}[h]
\centering
\includegraphics[width=0.95\linewidth]{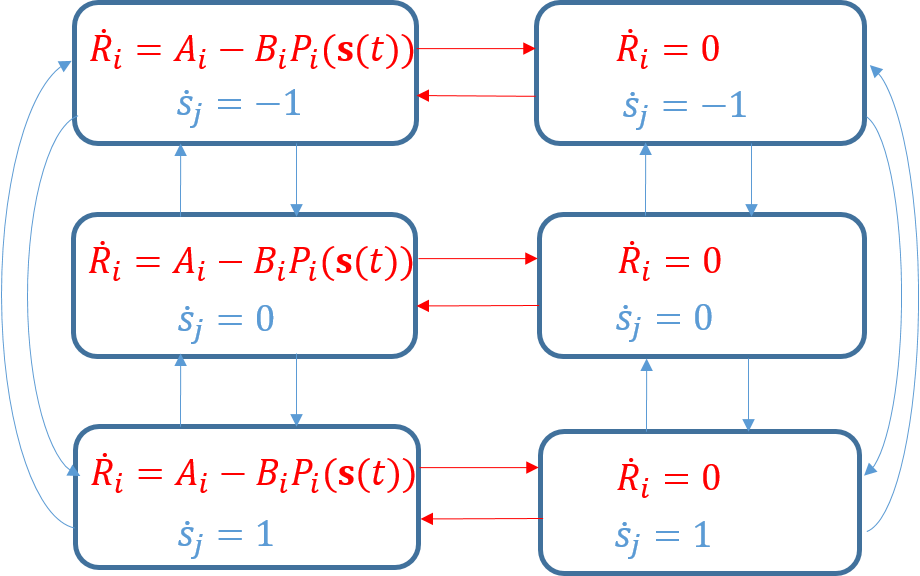}\caption{{\protect\small A
simple example of a hybrid system consisting of one agent and one target. The
system has six modes in which switches are triggered by events.}}%
\label{fig:Hybrid_sys}%
\end{figure}

A complete solution of the optimal control problem requires solving a TPBVP.
However, this is unnecessary, since the optimal control structure is fully
characterized by $u_{j}^{\ast}(t)\in\{1,0,-1\}$, it follows that we can
parameterize the optimal trajectory (illustrated in Fig. \ref{fig:para_illus})
so as to determine $(i)$ control switching points in $[0,L]$, where an agent
switches its control from $\pm1$ to $\mp1$ or possibly $0$ and $(ii)$
corresponding dwell times so that the cost in (\ref{eq:costfunction}) is
minimized. In other words, the optimal trajectory of each agent $j$ is fully
characterized by two parameter vectors: switching points $\bm\theta
_{j}=[\theta_{j1},\theta_{j2}...\theta_{j\Gamma}]$ and dwell times
$\mathbf{w}_{j}=[w_{j1},w_{j2}...w_{j\Gamma^{\prime}}]$ where $\Gamma$ and
$\Gamma^{\prime}$ depend on the given time horizon $T$. This defines a hybrid
system with state dynamics (\ref{eq:DynOfS}), (\ref{eq:multiDynR}). Figure
\ref{fig:Hybrid_sys} shows a simple example of such a hybrid system consisting
of one agent and one target. The dynamics remain unchanged in between events,
i.e., the points $\theta_{j1},\ldots,\theta_{j\Gamma}$ above and instants when
$R_{i}(t)$ switches from $>0$ to $0$ or vice versa. Therefore, the overall
cost function (\ref{eq:costfunction}) can be parametrically expressed as
$J(\bm\theta,\mathbf{w})$ and rewritten as the sum of costs over corresponding
inter-event intervals over a given time horizon:%
\begin{equation}
\min_{\mathbf{\bm\theta\in\Theta,\mathbf{w}\geq0}}\text{ \ }J(\bm\theta
,\mathbf{w})=\frac{1}{T}\sum_{k=0}^{K}\int_{\tau_{k}(\bm\theta,\mathbf{w}%
)}^{\tau_{k+1}(\bm\theta,\mathbf{w})}\sum_{i=1}^{M}R_{i}(t)dt
\label{eq:paramCost}%
\end{equation}
This allows us to apply IPA to determine a gradient $\nabla J(\bm\theta
,\mathbf{w})$ with respect to those parameters of the agent trajectories and
apply any standard gradient descent algorithm to obtain an optimal solution.

\section{Infinitesimal Perturbation Analysis\label{sec:IPA}}

We briefly review the IPA framework for general stochastic hybrid systems as
presented in \cite{cassandras2010perturbation}. Let $\{\tau_{k}(\theta)\}$,
$k=1,\ldots,K$, denote the occurrence times of all events in the state
trajectory of a hybrid system with dynamics $\dot{x}\ =\ f_{k}(x,\theta,t)$
over an interval $[\tau_{k}(\theta),\tau_{k+1}(\theta))$, where $\theta
\in\Theta$ is some parameter vector and $\Theta$ is a given compact, convex
set. For convenience, we set $\tau_{0}=0$ and $\tau_{K+1}=T$. We use the
Jacobian matrix notation: $x^{\prime}(t)\equiv\frac{\partial x(\theta
,t)}{\partial\theta}$ and $\tau_{k}^{\prime}\equiv\frac{\partial\tau
_{k}(\theta)}{\partial\theta}$, for all state and event time derivatives. It
is shown in \cite{cassandras2010perturbation} that
\begin{equation}
\frac{d}{dt}x^{\prime}(t)=\frac{\partial f_{k}(t)}{\partial x}x^{\prime
}(t)+\frac{\partial f_{k}(t)}{\partial\theta}, \label{eq:IPA_1}%
\end{equation}
for $t\in\lbrack\tau_{k},\tau_{k+1})$ with boundary condition:
\begin{equation}
x^{\prime}(\tau_{k}^{+})=x^{\prime}(\tau_{k}^{-})+[f_{k-1}(\tau_{k}^{-}%
)-f_{k}(\tau_{k}^{+})]\tau_{k}^{\prime} \label{eq:IPA_2}%
\end{equation}
for $k=1,...,K$. In order to complete the evaluation of $x^{\prime}(\tau
_{k}^{+})$ in (\ref{eq:IPA_2}), we need to determine $\tau_{k}^{\prime}$. If
the event at $\tau_{k}$ is \emph{exogenous} (i.e., independent of $\theta$),
$\tau_{k}^{\prime}=0$. However, if the event is \emph{endogenous}, there
exists a continuously differentiable function $g_{k}:\mathbb{R}^{n}%
\times\Theta\rightarrow\mathbb{R}$ such that $\tau_{k}\ =\ \min\{t>\tau
_{k-1}\ :\ g_{k}\left(  x\left(  \theta,t\right)  ,\theta\right)  =0\}$ and
\begin{equation}
\tau_{k}^{\prime}=-[\frac{\partial g_{k}}{\partial x}f_{k}(\tau_{k}^{-}%
)]^{-1}(\frac{\partial g_{k}}{\partial\theta}+\frac{\partial g_{k}}{\partial
x}x^{\prime}(\tau_{k}^{-})) \label{eq:IPA_3}%
\end{equation}
as long as $\frac{\partial g_{k}}{\partial x}f_{k}(\tau_{k}^{-})\neq0$
(details may be found in \cite{cassandras2010perturbation}).

In our setting, the cost along a given trajectory is $\frac{1}{T}\int_{0}%
^{T}\sum_{i=1}^{M}R_{i}(t)dt$. Following (\ref{eq:paramCost}), the gradient
for each agent $j$ denoted by $\nabla_{j}J(\bm\theta,\mathbf{w})=[\frac
{\partial J(\bm\theta,\mathbf{w})}{\partial\bm\theta_{j}},\frac{\partial
J(\bm\theta,\mathbf{w})}{\partial\mathbf{w}_{j}}]^{\text{T}}$ is
\begin{equation}
\nabla_{j}J(\bm\theta,\mathbf{w})=\frac{1}{T}\sum_{k=0}^{K}\sum_{i=1}^{M}%
\int_{\tau_{k}(\bm\theta,\mathbf{w})}^{\tau_{k+1}(\bm\theta,\mathbf{w})}%
\nabla_{j}R_{i}(t)dt \label{eq:GradientParametricObj}%
\end{equation}
where $\nabla_{j}R_{i}(t)=[\frac{\partial R_{i}(t)}{\partial\bm\theta_{j}%
},\frac{\partial R_{i}(t)}{\partial\mathbf{w}_{j}}]^{\text{T}}$.

We begin by deriving the gradient above within any inter-event interval
$[\tau_{k},\tau_{k+1})$ when the dynamics of both agent $j$ and target $i$
remain unchanged. Then, in Section \ref{sec:Events}, we will define all events
involved in switching these dynamics, hence affecting the gradient
evaluation possibly through discontinuities characterized by (\ref{eq:IPA_2}).
We proceed with the derivation of $\frac{\partial R_{i}(t)}{\partial
\bm\theta_{j}}$, since $\frac{\partial R_{i}(t)}{\partial\mathbf{w}_{j}}$ can
be derived in a similar way.

It follows from \eqref{eq:IPA_1}, observing that the first term vanishes since
$f_{k}(t)=\dot{R}_{i}(t)$ is not an explicit function of $R_{i}(t)$, that
$
\frac{d}{dt}\frac{\partial R_{i}(t)}{\partial\bm\theta_{j}}=\frac{\partial
\dot{R}_{i}(t)}{\partial\bm\theta_{j}}%
$.
Then, in view of (\ref{eq:multiDynR}), we have for all $t\in\lbrack\tau
_{k},\tau_{k+1})$:
\begin{equation}
\frac{\partial R_{i}(t)}{\partial\bm\theta_{j}}=\frac{\partial R_{i}(\tau
_{k}^{+})}{\partial\bm\theta_{j}}\text{ \ \ if }\dot{R}_{i}(t)=0
\label{eq:IPA_dotR_0}%
\end{equation}
and
\begin{equation}
\begin{split}
\frac{\partial R_{i}(t)}{\partial\bm\theta_{j}}=\frac{\partial R_{i}(\tau
_{k}^{+})}{\partial\bm\theta_{j}}-B_{i}\int_{\tau_{k}}^{t}\frac{\partial
P_{i}(\mathbf{s}(\tau))}{\partial\theta_{j}}d\tau\text{ \ \ }%
\label{eq:IPA_dotR_ne_0}\\
\text{ if }\dot{R}_{i}(t)=A_{i}-B_{i}P_{i}(\mathbf{s}(t))
\end{split}
\end{equation}

The integrand in \eqref{eq:IPA_dotR_ne_0} is obtained from
(\ref{JointSensingProb}):
\begin{equation}
\frac{\partial P_{i}(\mathbf{s}(\tau))}{\partial\theta_{j}}=\frac{\partial
p_{ij}(s_{j}(\tau))}{\partial s_{j}}\frac{\partial s_{j}(\tau)}{\partial
\bm\theta_{j}}\prod_{\substack{{g}\in\mathcal{B}_{i}(\tau)\\{g\neq j}}}\left[
1-p_{ig}\left(  s_{g}(\tau)\right)  \right]  \label{eq:diff_Pi}%
\end{equation}
Note that $\frac{\partial p_{ij}\left(  s_{j}(\tau)\right)  }{\partial s_{j}}$
is piece-wise constant and takes values in $\{0,\pm\frac{1}{r_{j}}\}$
depending on $|s_{j}(t)-x_{i}|$ and $r_{j}$ (see agent sensing mode \eqref{eq:sensing_single}). We can, therefore, factor the constant
$\frac{\partial p_{ij}\left(  s_{j}(\tau)\right)  }{\partial s_{j}}$ out of
the integral in \eqref{eq:IPA_dotR_ne_0}. As for the term $\frac{\partial
s_{j}(\tau)}{\partial\bm\theta_{j}}$, we apply \eqref{eq:IPA_1} and
(\ref{eq:DynOfS}) to obtain
\begin{equation}
\frac{d}{dt}\frac{\partial s_{j}(\tau)}{\partial\bm\theta_{j}}=0
\label{eq:ps_const}%
\end{equation}
Therefore, $\frac{\partial s_{j}(\tau)}{\partial\bm\theta_{j}}=\frac{\partial
s_{j}(\tau_{k}^{+})}{\partial\bm\theta_{j}}$ which is also a constant. The
product term in \eqref{eq:diff_Pi} captures the contributions from all agents
other than $j$ in monitoring target $i$. Using the definition of
$\mathcal{N}_{ij}(t)$ in \eqref{eq:new_def_neighbor}, it can be restricted to
this set, since for any agent $g\not \in \mathcal{N}_{ij}(t)$, $p_{ig}%
(s_{g}(t))=0$. For notational simplicity, we define the integral of this term
over $[\tau_{k},t)$, $t<\tau_{k+1}$, as:
\begin{equation}
G_{ij}(t)=\int_{\tau_{k}}^{t}\prod_{g\in\mathcal{N}_{ij}(\tau)}\left[
1-p_{ig}\left(  s_{g}(\tau)\right)  \right]  d\tau\label{eq:collaborationTerm}%
\end{equation}
which can be interpreted as a \textquotedblleft collaboration
factor\textquotedblright\ involving all agents in $\mathcal{N}_{ij}(\tau)$.
Clearly, this is affected by an agent leaving or joining the neighbor set
$\mathcal{N}_{ij}(\tau)$ which motivates defining an event associated with
such changes (see Section \ref{sec:Events}).

When we combine \eqref{eq:IPA_dotR_0} and \eqref{eq:IPA_dotR_ne_0}, the
derivative $\frac{\partial R_{i}(t)}{\partial\bm\theta_{j}}$, $i=1,\ldots,M$,
over any inter-event interval $[\tau_{k},\tau_{k+1})$ becomes:
\begin{equation}
\frac{\partial R_{i}(t)}{\partial\bm\theta_{j}}\hspace{-1mm}=\hspace
{-1mm}\frac{\partial R_{i}(\tau_{k}^{+})}{\partial\bm\theta_{j}}-%
\begin{cases}
0\qquad\text{if }R_{i}(t)=0,\text{ }A_{i}\leq B_{i}P_{i}(\mathbf{s}(t))\\
\!\!B_{i}\frac{\partial p_{ij}(\!s_{j}(\tau_{k}^{+})\!)}{\partial s_{j}}\frac{\partial
s_{j}(\tau_{k}^{+})}{\partial\bm\theta_{j}}G_{ij}(t)\text{ otherwise}%
\end{cases}
\label{eq:pRptheta}%
\end{equation}
A similar derivation can be applied to the derivative $\frac{\partial
R_{i}(t)}{\partial\mathbf{w}_{j}}$ and gives:
\begin{equation}
\frac{\partial R_{i}(t)}{\partial\mathbf{w}_{j}}\hspace*{-1mm}=\hspace*{-1mm}\frac{\partial R_{i}(\tau
_{k}^{+})}{\partial\mathbf{w}_{j}}-%
\begin{cases}
0\qquad\text{if }R_{i}(t)=0,\text{ }A_{i}\leq B_{i}P_{i}(\mathbf{s}(t))\\
\!\!B_{i}\frac{\partial p_{ij}(\!s_{j}(\tau_{k}^{+})\!)}{\partial s_{j}}\frac{\partial s_{j}(\tau_{k}^{+})}{\partial\mathbf{w}_{j}}G_{ij}(t)\text{ otherwise}%
\end{cases}
\label{eq:pRpw}%
\end{equation}

\subsection{Events in the hybrid system\label{sec:Events}}

We are now in a position to define as \textquotedblleft
events\textquotedblright\ all switches in the hybrid system which can result
in changes in the derivatives in \eqref{eq:pRptheta} and \eqref{eq:pRpw} so we
can apply (\ref{eq:IPA_2}) to determine the initial conditions $\frac{\partial
R_{i}(\tau_{k}^{+})}{\partial\bm\theta_{j}}$ and $\frac{\partial R_{i}%
(\tau_{k}^{+})}{\partial\mathbf{w}_{j}}$ at $t=\tau_{k}$, as well as the terms
$\frac{\partial s_{j}(\tau_{k}^{+})}{\partial\bm\theta_{j}}$ and
$\frac{\partial s_{j}(\tau_{k}^{+})}{\partial\mathbf{w}_{j}}$.


We classify events into four categories depending on the effect they have on
target dynamics (type I), agent sensing relative to a target (type II), agent
dynamics (type III), and neighbor set $N_{ij}(t)$ (type IV). Referring to Fig.
\ref{fig:Hybrid_sys}), observe that only event types I and III (red and blue
arrows) affect the dynamics of the corresponding target and agent. Event types
II and IV do not change the system dynamics but still affect the derivative
values in \eqref{eq:pRptheta} and \eqref{eq:pRpw}. In what follows, we define
all events types and their corresponding effects on \eqref{eq:pRptheta} and
\eqref{eq:pRpw} and summarize them in Table \ref{table:events}.

\begin{table}[ptb]
\caption{{\protect\small Events in agent-target system}}%
\label{table:events}%
\centering
\begin{tabular}
[c]{|p{1.5cm}||p{5.7cm}|}\hline
Event Name & Description\\\hline
$\rho_{i}^{0}$ & $R_{i}(t)$ hits 0\\\hline
$\rho_{i}^{+}$ & $R_{i}(t)$ leaves 0\\\hline
$\pi_{ij}^{0}$ & $p_{ij}(s_{j}(t))$ hits 0\\\hline
$\pi_{ij}^{+}$ & $p_{ij}(s_{j}(t))$ leaves 0\\\hline
$\nu_{j}^{(1,0)}$ & $u_{j}(t)$ switches from $1$ to $0$\\\hline
$\nu_{j}^{(-1,0)}$ & $u_{j}(t)$ switches from $-1$ to $0$\\\hline
$\nu_{j}^{(0,1)}$ & $u_{j}(t)$ switches from $0$ to $1$\\\hline
$\nu_{j}^{(0,-1)}$ & $u_{j}(t)$ switches from $0$ to $-1$\\\hline
$\nu_{j}^{(1,-1)}$ & $u_{j}(t)$ switches from $1$ to $-1$\\\hline
$\nu_{j}^{(-1,1)}$ & $u_{j}(t)$ switches from $-1$ to $1$\\\hline
$\Delta_{ij}^{+}$ & $\mathcal{N}_{ij}(\tau^{+}) = \mathcal{N}_{ij}(\tau^{-})
\cup\{k \}$, $k \not \in \mathcal{N}_{ij}(\tau^{-})$\\\hline
$\Delta_{ij}^{-}$ & $\mathcal{N}_{ij}(\tau^{+}) = \mathcal{N}_{ij}(\tau^{-})
\setminus\{k \}$, $k \in\mathcal{N}_{ij}(\tau^{-})$\\\hline
\end{tabular}
\raggedright{Note: events in the table include all $i = 1, \ldots, M$ and $j =
1 ,\ldots, N$}
\end{table}

\textbf{Event type I: switches in target dynamics $\dot{R}_{i}(t)$.} Referring
to (\ref{eq:multiDynR}), when $R_{i}(t)$ either reaches zero or leaves zero,
the IPA derivative switches between \eqref{eq:IPA_dotR_0} and
\eqref{eq:IPA_dotR_ne_0}. We denote the former event as $\rho_{i}^{0}$ and the
latter as $\rho_{i}^{+}$ for all $i=1,\ldots,M$ (see Table \ref{table:events}%
). When such events occur, the dynamics of $s_{j}(t)$ in \eqref{eq:DynOfS}
remain unchanged, so it follows from (\ref{eq:IPA_2}) that $\nabla_{j}%
s_{j}(\tau_{k}^{-})=\nabla_{j}s_{j}(\tau_{k}^{+})$. However, the target
dynamics switch between $\dot{R}_{i}=A_{i}-B_{i}P_{i}(\mathbf{s}(t))$ and
$\dot{R}_{i}=0$ and cause discontinuities in $\nabla_{j}R_{i}(t)$ as follows.

\emph{Event $\rho_{i}^{0}$}: This event causes a transition from $\dot{R}%
_{i}(t)=A_{i}-B_{i}P_{i}(\mathbf{s}(t)),$ $t<\tau_{k}$ to $\dot{R}_{i}(t)=0,$
$t\geq\tau_{k}$. It is an endogenous event because its occurrence depends on
the parameters $\bm\theta,\mathbf{w}$ which dictate switches in $\mathbf{s}%
(t)$. We first evaluate $\tau_{k}^{\prime}$ from \eqref{eq:IPA_3} with
$g_{k}(R_{i}(t),t)=R_{i}(t)=0$ to get
\begin{equation}
\tau_{k}^{\prime}=-\frac{\nabla_{j}R_{i}(\tau_{k}^{-})}{A_{i}-B_{i}%
P_{i}(\mathbf{s}(\tau_{k}^{-}))} \label{eq:event1_1}%
\end{equation}
and then apply \eqref{eq:IPA_2} to obtain
\begin{equation}
\nabla_{j}R_{i}(\tau_{k}^{+})=\nabla_{j}R_{i}(\tau_{k}^{-})+\left[
A_{i}-B_{i}P_{i}(\mathbf{s}(\tau_{k}^{-}))-0\right]  \tau_{k}^{\prime}
\label{eq:event1_2}%
\end{equation}
Combining \eqref{eq:event1_1} and \eqref{eq:event1_2}, we get
\begin{equation}
\nabla_{j}R_{i}(\tau_{k}^{+})=0\quad\text{if event $\rho_{i}^{0}$ occurs at
$\tau_{k}$} \label{eq:tar_reset}%
\end{equation}

\emph{ Event $\rho_{i}^{+}$:} This event causes a transition from $\dot{R}%
_{i}(t)=0,$ $t<\tau_{k}$ to $\dot{R}_{i}(t)=A_{i}-B_{i}P_{i}(\mathbf{s}(t)),$
$t\geq\tau_{k}$. It is easy to see that the dynamics in both
(\ref{eq:DynOfS}) and (\ref{eq:multiDynR}) are continuous when this happens
and since $A_{i}-B_{i}P_{i}(\mathbf{s}(\tau_{k}))=0$ we have $\dot{R}_{i}%
(\tau_{k}^{-})=\dot{R}_{i}(\tau_{k}^{+})=0$. It follows from \eqref{eq:IPA_2}
that $\nabla_{j}R_{i}(\tau_{k}^{+})=\nabla_{j}R_{i}(\tau_{k}^{-})$. Moreover,
since $R_{i}(t)=0$, $\dot{R}_{i}(t)=0,$ $t<\tau_{k}$, we have $\nabla_{j}%
R_{i}(\tau_{k}^{-})=0$ and we get
\begin{equation}
\nabla_{j}R_{i}(\tau_{k}^{+})=0\quad\text{if event $\rho_{i}^{+}$ happens at
$\tau_{k}$} \label{eq:tar_plus}%
\end{equation}

\textbf{Remark 1:} Combining (\ref{eq:tar_reset}) and (\ref{eq:tar_plus}) with
\eqref{eq:pRptheta} and \eqref{eq:pRpw}, we conclude that a $\rho_{i}^{0}$
event occurring at $t=\tau_{k}$ resets the value of $\nabla_{j}R_{i}(t)$ to
$\nabla_{j}R_{i}(t)=0$ for all $j=1,\ldots,N$ regardless of the value
$\nabla_{j}R_{i}(\tau_{k}^{-})$ and the state of the agents. Moreover,
$R_{i}(t)=0$ and $\nabla_{j}R_{i}(t)=0$ for $t>\tau_{k}$ until the next
$\rho_{i}^{+}$ event occurs.

\textbf{Event type II: switches in agent sensing $p_{ij}(s_{j}(t))$.} These
events trigger a switch in $\frac{\partial p_{ij}(s_{j}(t))}{\partial s_{j}}$
from $\pm\frac{1}{r_{j}}$ to $0$ or vice versa in \eqref{eq:pRptheta} and
\eqref{eq:pRpw}. We denote the former event as $\pi_{ij}^{0}$ and the latter
as $\pi_{ij}^{+}$. These events trigger a switch of $\frac{\partial
p_{ij}(s_{j}(t))}{\partial s_{j}}$ at $t=\tau_{k}$ from $\pm\frac{1}{r_{j}}$
to $0$ or vice versa in \eqref{eq:pRptheta} and \eqref{eq:pRpw}. However, the
dynamics in both (\ref{eq:DynOfS}) and (\ref{eq:multiDynR}) remain unchanged
when this happens (due to the continuity of the sensing function
$p_{ij}\left(  s_{j}(t)\right)  $) and it follows from (\ref{eq:IPA_2}) that
$\nabla_{j}R_{i}(\tau_{k}^{+})=\nabla_{j}R_{i}(\tau_{k}^{-})\text{ and }%
\nabla_{j}s_{j}(\tau_{k}^{+})=\nabla_{j}s_{j}(\tau_{k}^{-})$.

\textbf{Event type III: switches in agent dynamics $\dot{s}_{j}(t)$.}
Referring to (\ref{eq:DynOfS}), these are events that cause a switch in the
optimal control values $u_{j}^{\ast}(\tau_{k})$: $(i)$ $\pm1\rightarrow0$,
$(ii)$ $0\rightarrow\pm1$, and $(iii)$ $\pm1\rightarrow\mp1$. We denote these
events as $\nu_{j}^{(1,0)},\nu_{j}^{(-1,0)},\nu_{j}^{(0,1)},$ $\nu
_{j}^{(0,-1)},\nu_{j}^{(-1,1)},\nu_{j}^{(1,-1)}$ using the general notation
$\nu_{j}^{(\ast,\ast)}$ with the superscript corresponding to the six total
possible control switches. The effect of these events in (\ref{eq:pRptheta})
and (\ref{eq:pRpw}) is through possible discontinuities in the terms
$\frac{\partial s_{j}(t)}{\partial\bm\theta_{j}}$ and $\frac{\partial
s_{j}(t)}{\partial\mathbf{w}_{j}}$ at $t=\tau_{k}$. Clearly, the gradient
cannot be affected by future events, so we consider all prior and current
control switches indexed by $l=1,2...,\xi$ where $\xi$ is the current control
switch and $\theta_{jl}$, $w_{jl}$ are the $l$-th switching point and dwelling
time respectively. These agent control switches are endogenous events with
switching functions $g_{k}(s_{j}(t),t)=s_{j}-\theta_{jl}=0$. We can now apply
(\ref{eq:IPA_2}) and (\ref{eq:IPA_3}) to \eqref{eq:DynOfS}, similar to the
derivation for type I events. We omit the details (which can be found in
\cite{cassandras2013optimal}) and present the final results.

\emph{Events }$\nu_{j}^{(1,0)},\nu_{j}^{(-1,0)}$: These are switches such that
$u_{j}(\tau_{k}^{-})=\pm1$, $u_{j}(\tau_{k}^{+})=0$ and we get
\begin{align}
&  \frac{\partial s_{j}}{\partial\theta_{jl}}(\tau_{k}^{+})=%
\begin{cases}
1 & \text{if }l=\xi\\
0 & \text{if }l<\xi
\end{cases}
\label{eq:Event2_11}\\
&  \frac{\partial s_{j}}{\partial\omega_{jl}}(\tau_{k}^{+})=0\quad\text{for
all }l\leq\xi\label{eq:Event2_12}%
\end{align}

\emph{Events }$\nu_{j}^{(0,1)},$ $\nu_{j}^{(0,-1)}$: These are switches such
that $u_{j}(\tau_{k}^{-})=0,$ $u_{j}(\tau_{k}^{+})=\pm1$ and we get
\begin{equation}
\frac{\partial s_{j}}{\partial\theta_{jl}}(\tau_{k}^{+})\hspace{-1mm}%
=\hspace{-1mm}%
\begin{cases}
\frac{\partial s_{j}}{\partial\theta_{jl}}(\tau_{k}^{-})-\hspace{-1mm}%
u_{j}(\tau_{k}^{+})sgn\big(\theta_{j\xi}-\theta_{j(\xi-1)}\big)\text{ if
}l=\xi\\
\frac{\partial s_{j}}{\partial\theta_{jl}}(\tau_{k}^{-})-u_{j}(\tau_{k}%
^{+})\Big[sgn(\theta_{jl}-\theta_{j(l-1)})\\
\qquad\qquad-sgn(\theta_{j(l+1)}-\theta_{jl})\Big]\qquad\quad\text{if }l<\xi
\end{cases}
\label{eq:Event2_21}%
\end{equation}%
\begin{equation}
\frac{\partial s_{j}}{\partial w_{jl}}(\tau_{k}^{+})=-u_{j}(\tau_{k}^{+}%
)\quad\text{for all }l\leq\xi\label{eq:Event2_22}%
\end{equation}

\emph{Events }$\nu_{j}^{(-1,1)},\nu_{j}^{(1,-1)}$: These are switches such
that $u_{j}(\tau_{k}^{-})=\pm1,$ $u_{j}(\tau_{k}^{+})=\mp1$ so that a dwell
time is not involved and we get
\begin{equation}
\frac{\partial s_{j}}{\partial\theta_{jl}}(\tau_{k}^{+})=%
\begin{cases}
2 & \text{if }l=\xi\\
-\frac{\partial s_{j}}{\partial\theta_{jl}}(\tau_{k}^{-}) & \text{if }l<\xi
\end{cases}
\label{eq:Event2_3}%
\end{equation}

\textbf{Remark 2: }Observe that $\nabla_{j}s_{j}(t)$ is independent of the
states of other agents $k\neq j$. This follows from the fact that $\nabla
_{j}s_{j}(t)$ is constant over inter-event intervals $[\tau_{k},\tau_{k+1})$
as shown in \eqref{eq:ps_const} and only depends on parameter and control
values known to agent $j$ as seen in \eqref{eq:Event2_11} -
\eqref{eq:Event2_3}. Moreover, if $k\neq j$, $\nabla_{k}s_{j}(t)=0$.

\textbf{Event type IV: changes in neighbor sets $\mathcal{N}_{ij}(t)$.} These
events change the topology of the agent-target network by altering the
neighbors of agent $j$, hence affecting the value of $G_{ij}(t)$ in
\eqref{eq:collaborationTerm} which in turn affects (\ref{eq:pRptheta}) and
(\ref{eq:pRpw}). We denote by $\Delta_{ij}^{+}$ the event causing the addition
of an agent to the neighbor set $\mathcal{N}_{ij}(t)$ and by $\Delta_{ij}^{-}$
the event causing the removal of an agent from the neighbor set $\mathcal{N}%
_{ij}(t)$. However, the dynamics of both $R_{i}(t)$ and $s_{j}(t)$ remain
unchanged when these events occur. Due to the continuity of the sensing
function $p_{ig}\left(  s_{g}(\tau)\right)  $ in \eqref{eq:collaborationTerm},
the addition/removal of an agent $g$ to/from the set $\mathcal{N}_{ij}(\tau)$
does not affect the continuity of $G_{ij}(t)$, which implies $\nabla_{j}%
R_{i}(\tau_{k}^{+})=\nabla_{j}R_{i}(\tau_{k}^{-})$ as well as $\nabla_{j}%
s_{j}(\tau_{k}^{+})=\nabla_{j}s_{j}(\tau_{k}^{-})$.

The set of all events defined above and summarized in Table \ref{table:events}
is denoted by $\mathcal{E}$. Furthermore, we define the set of all type III
events of the form $\nu_{j}^{(\ast,\ast)}$ as the \emph{agent event set}
$\mathcal{E}^{A}$ and the set of all other events (type I, III, and IV) as the
\emph{target event set} $\mathcal{E}^{T}$. The subset of $\mathcal{E}^{A}$
that contains only events related to agent $j$ is denoted by $\mathcal{E}%
_{j}^{A}$. Similarly, the subset of $\mathcal{E}^{T}$ that contains only
events related to target $i$ is denoted by $\mathcal{E}_{i}^{T}$. We then have:

\begin{definition}
The local event set of any agent $j$ is the union of agent events
$\mathcal{E}_{j}^{A}$ and target events $\mathcal{E}_{i}^{T}$ for all
$i\in\mathcal{T}_{j}(t)$:
\begin{equation}
\mathcal{E}_{j}(t)=\mathcal{E}_{j}^{A} \bigcup_{i\in\mathcal{T}_{j}(t)}
\mathcal{E}_{i}^{T} \label{LocalEventSet}%
\end{equation}

\end{definition}

In contrast, the global event set for agent $j$ includes all non-neighboring
target events in $\mathcal{E}_{i}^{T}\text{ for all }i\not \in \mathcal{T}%
_{j}$ and non-neighboring agent events $\mathcal{E}_{k}^{A},\text{ for all
}k\not \in \mathcal{A}_{j}$. Based on the limited information model of Section
\ref{sec:PM_form}, we define the \emph{local information set} of agent $j$,
denoted by $I_{j}(t)$, as follows:

\begin{definition}
The local information set of any agent $j$ is the union of its local event set
and those of its neighbors in $\mathcal{N}_{ij}(t)$ for all $i\in
\mathcal{T}_{j}(t)$:
\begin{equation}
I_{j}(t)=\mathcal{E}_{j}(t)\bigcup_{k\in\mathcal{N}_{ij}(t),i\in
\mathcal{T}_{j}(t)}\mathcal{E}_{k}(t). \label{LocalInfoSet}%
\end{equation}

\end{definition}

This includes all local information necessary for agent $j$ to evaluate the
IPA gradient $\nabla_{j}R_{i}(t)$ for $i\in\mathcal{T}_{j}(t)$. Observe that
agent $j$ does not need to communicate with all its neighbors in
$\mathcal{A}_{j}(t)$, but only a subset which includes those neighbors who are
sharing the same target(s) as $j$ at time $t$ since $\bigcup_{i\in
\mathcal{T}_{j}(t)}\mathcal{N}_{ij}(t)\subseteq\mathcal{A}_{j}(t)$.

\textbf{Remark 3: }It is clear from the analysis thus far, that IPA is
entirely \emph{event-driven}, since all gradient updates happen exclusively at
events occurring at times $\tau_{k}({\bm\theta},\mathbf{w})$, $k=1,2,\ldots$.
Thus, this approach scales with the number of events characterizing the hybrid
system, and not its (generally much larger) state space.

\section{Event-Driven decentralized gradient evaluation and optimization}
\label{sec:DecentOpt}
Our main results are presented in this section. In particular, we show in
Theorem 1 that each agent can evaluate the gradient of the objective function
in (\ref{eq:paramCost}) with respect to its own controllable parameters
$\bm\theta_{j}$ and $\mathbf{w}_{j}$ based on its local information set
(\ref{LocalInfoSet}) and only one non-local event. We begin with the following
lemma which asserts that the gradient $\nabla_{j}R_{i}(t)$ takes a very simple
form as long as $i\notin\mathcal{T}_{j}(t)$, i.e., while target $i$ cannot be
sensed by agent $j$.

\begin{lemma}
\label{lemma1} Let $t\in\lbrack t_{1},t_{2}]$ such that $i\not \in
\mathcal{T}_{j}(t)$. Then,

\begin{enumerate}
\item If $R_{i}(t)>0$ for all $t\in[t_{1},t_{2}]$, then
\begin{equation}
\nabla_{j}R_{i}(t)=\nabla_{j}R_{i}(t_{1}^{+}) \label{eq:lemma1_1}%
\end{equation}

\item If there exists an event $\rho_{i}^{0}$ at $\tau\in(t_{1},t_{2})$, then
\begin{equation}
\nabla_{j}R_{i}(t)=%
\begin{cases}
\nabla_{j}R_{i}(t_{1}^{+}) & t\in [t_{1},\tau)\\
0 & t\in\lbrack\tau,t_{2}]
\end{cases}
\label{eq:lemma1_2}%
\end{equation}

\end{enumerate}
\end{lemma}

\emph{Proof:} By the definition of $\mathcal{T}_{j}(t)$, when $i\not \in
\mathcal{T}_{j}(t)$ we have $\Vert s_{j}(t)-x_{i}\Vert>r_{j}$ and
$\frac{\partial p_{ij}(s_{j}(t))}{\partial s_{j}}=0$ for all $t\in\lbrack
t_{1},t_{2}]$. If $R_{i}(t)>0$ for all $t\in\lbrack t_{1},t_{2}]$, it follows
directly from \eqref{eq:pRptheta} and \eqref{eq:pRpw} that $\nabla_{j}%
R_{i}(t)=\nabla_{j}R_{i}(t_{1}^{+})$. Otherwise, there exists an event
$\rho_{i}^{0}$ at time $\tau\in(t_{1},t_{2})$ which results in $R_{i}(\tau
)=0$. The previous argument applies to $(t_{1},\tau)$ giving $\nabla_{j}%
R_{i}(t)=\nabla_{j}R_{i}(t_{1}^{+})$ for $t\in\lbrack t_{1},\tau)$. According
to \eqref{eq:tar_reset}, event $\rho_{i}^{0}$ resets the gradient to
$\nabla_{j}R_{i}(\tau)=0$. Subsequently, over $[\tau,t_{2}]$, regardless of
which of the cases in \eqref{eq:pRptheta} and \eqref{eq:pRpw} applies, it
holds that $\nabla_{j}R_{i}(t)=0$. $\blacksquare$

\begin{corollary}
$\nabla_{j} R_{i}(t)$ is independent of events $\rho_{i}^{+}$ for
$i\not \in \mathcal{T}_{j}(t)$.
\end{corollary}

\emph{Proof:} Note that the $\rho_i^+$ event can only occur after a $\rho_i^0$ event. The proof is self-evident following Lemma \ref{lemma1}. We have
$\nabla_{j}R_{i}(t)=0$ for $t> \tau$ until target $i$ joins the target neighborhood of agent
$j$. Therefore, any non-local $\rho_{i}^{+}$ event that may occur cannot
affect $\nabla_{j}R_{i}(t)$. $\blacksquare$

Lemma \ref{lemma1} and its Corollary imply that agent $j$ does not need any
knowledge of non-neighboring target events except for $\rho_{i}^{0}$ with
$i\not \in \mathcal{T}_{j}(t)$ in order to evaluate its gradient. We can
further establish that the gradient $\nabla_{j}J(\bm\theta,\mathbf{w})$ along
the agent trajectory is affected by only local events in $I_{j}(t)$, as
defined in (\ref{LocalInfoSet}), and a small subset of global events.

\begin{lemma}
\label{lemma2} A sufficient event set to evaluate $\nabla_{j}J(\bm\theta
,\mathbf{w})$ is $I_{j}(t)\cup\{\rho_{i}^{0}:i\not \in \mathcal{T}_{j}(t)\}$.
\end{lemma}

\emph{Proof:} Let $\tau_{k}$ be any event time when $\mathcal{T}_{j}(\tau
_{k})$ is altered, i.e., a new target is added to the target neighborhood of
agent $j$ or one is removed from it. From Lemma \ref{lemma1}, if
$i\not \in \mathcal{T}_{j}(t)$, then either $\nabla_{j}R_{i}(t)=\nabla
_{j}R_{i}(\tau_{k})$ and remains constant at this value or $\nabla_{j}%
R_{i}(t)=0$, depending on whether an event $\rho_{i}^{0}$ takes place. It
follows from (\ref{eq:GradientParametricObj}) that the objective function
gradient can be rewritten as
\begin{align}
\nabla_{j}J(\bm\theta,\mathbf{w})  &  =\sum_{k=0}^{K}\sum_{i=1}^{M}\int
_{\tau_{k}}^{\tau_{k+1}}\nabla_{j}R_{i}(t)dt\nonumber\\
&  =\sum_{k=0}^{K}\Big(\sum_{i\not \in \mathcal{T}_{j}(\tau_{k})}\nabla
_{j}R_{i}(\tau_{k})(\tau_{k+1}-\tau_{k})\nonumber\\
&  \qquad\qquad+\sum_{i\in\mathcal{T}_{j}(\tau_{k})}\int_{\tau_{k}}%
^{\tau_{k+1}}\nabla_{j}R_{i}(t)dt\Big) \label{eq:lemma2}%
\end{align}
The value of $\nabla_{j}R_{i}(\tau_{k})$ in the first term of
\eqref{eq:lemma2} depends on $\{\rho_{i}^{0}:i\not \in \mathcal{T}_{j}(t)\}$
which is a subset of events non-local to agent $j$. The second term of
\eqref{eq:lemma2} depends only on the local information set events $I_{j}(t)$
since target $i\in\mathcal{T}_{j}(t)$ is local to agent $j$. Therefore,
$I_{j}(t)\cup\{\rho_{i}^{0}:i\not \in \mathcal{T}_{j}(t)\}$ is a sufficient
event set to evaluate $\nabla_{j}J(\bm\theta,\mathbf{w})$. $\blacksquare$

\textbf{Remark 4:} Although an event $\rho_{i}^{0}$ for $i\not \in
\mathcal{T}_{j}(t)$ is non-local to agent $j$, it must be observed by at least
one agent $k\neq j$ such that $i\in\mathcal{T}_{k}(t)$. This is because
$\rho_{i}^{0}$ at some time $\tau_{k}$ can only take place if one or more
agents in its neighborhood cause a transition from $R_{i}(\tau_{k}^{-})>0$ to
$R_{i}(\tau_{k})=0$ in (\ref{eq:multiDynR}). Therefore, such events can be
communicated to agent $j$ through the agent network, possibly with some
delays. The implication of Lemma \ref{lemma2} is an \textquotedblleft almost
decentralized\textquotedblright\ algorithm in which each agent optimizes its
trajectory through the gradient $\nabla_{j}J(\bm\theta,\mathbf{w})$ using only
agent local information; the only exception is occasional target uncertainty
depletion events transmitted to it from other agents.

Returning to the parametric optimization problem (\ref{eq:paramCost}), a
centralized solution was obtained in \cite{zhou2016optimal} using the IPA
gradients in \eqref{eq:pRptheta} and \eqref{eq:pRpw} and a standard gradient
descent scheme to optimize the parameter vector $[\bm\theta,\mathbf{w}%
]^{\text{T}}$ as follows:%
\begin{equation}
\left[  \bm\theta^{l+1},\mathbf{w}^{l+1}\right]  ^{\text{T}}=\left[
\bm\theta^{l},\mathbf{w}^{l}\right]  ^{\text{T}}-[\alpha_{\theta}^{l}%
,\alpha_{w}^{l}]\nabla J(\bm\theta,\mathbf{w}) \label{CentralGradAlgo}%
\end{equation}
where $l=0,1,\ldots$ is the iteration index and $\alpha_{\theta}^{l}$ and
$\alpha_{w}^{l}$ are diminishing step-size sequences satisfying $\sum
_{l=0}^{\infty}\alpha_{\theta}^{l}=\infty,\lim_{l\rightarrow\infty}%
\alpha_{\theta}^{l}=0$ and $\sum_{l=0}^{\infty}\alpha_{w}^{l}=\infty
,\lim_{l\rightarrow\infty}\alpha_{w}^{l}=0$. A decentralized version of
(\ref{CentralGradAlgo}) is
\begin{equation}
\left[  \bm\theta_{j}^{l+1},\mathbf{w}_{j}^{l+1}\right]  ^{\text{T}}=\left[
\bm\theta_{j}^{l},\mathbf{w}_{j}^{l}\right]  ^{\text{T}}-[\alpha_{\theta}%
^{l},\alpha_{w}^{l}]\nabla_{j}J(\bar{\bm\theta},\bar{\mathbf{w}})
\label{eq:ParaUpdate}%
\end{equation}
where $\bar{\bm\theta}$ and $\bar{\mathbf{w}}$ are agent $j$'s estimates based
on the limited information provided in Lemma \ref{lemma2}.

\begin{theorem}
\label{thm1} Any centralized solution of (\ref{eq:paramCost}) through
(\ref{CentralGradAlgo}) can be recovered by (\ref{eq:ParaUpdate}) in which
each agent $j$ optimizes its trajectory given the following conditions:

$1)$ Initial parameters $[\bm\theta_{j}^{0},\mathbf{w}_{j}^{0}]$;

$2)$ The local information set $I_{j}(t)$;

$3)$ The subset of the global information set $\{\rho_{i}^{0},i\not \in
\mathcal{T}_{j}(t)\}$.
\end{theorem}

\emph{Proof:} The proof is immediate from Lemma \ref{lemma2}. The gradient
$\nabla_{j}J(\bm\theta,\mathbf{w})$ can be evaluated by each agent given
conditions $2$ and $3$. Condition $1$ provides initial parameters for each
agent trajectory in order to execute (\ref{eq:ParaUpdate}). $\blacksquare$

Note that condition $3$ involves only a small subset of global events. As
shown in our simulation results in Section \ref{Simulation}, ignoring such
non-local events will affect the cooperation among agents and increase the
final cost. Thus, it can be interpreted as the \textquotedblleft price of
anarchy\textquotedblright\ commonly associated with decentralization limiting
agent actions to only local information.

It is important to point out that the method of Theorem 1 relies on the
gradient $\nabla_{j}R_{i}(t)$ for $i\not \in \mathcal{T}_{j}(t)$ and not on
$R_{i}(t)$. In fact, there is no attempt by agent $j$ to reconstruct or
estimate the states of targets $i\not \in \mathcal{T}_{j}(t)$; the only
information from such targets is provided through the occasional $\rho_{i}%
^{0}$ events.

We briefly discuss next some open issues defining ongoing research directions. While the event-driven nature of IPA has several computational advantages (see \textbf{Remark 3}), the optimization process depends on these events being observed so as to \textquotedblleft excite\textquotedblright\ algorithms such as
(\ref{CentralGradAlgo}) and (\ref{eq:paramCost}). To resolve this \emph{event
excitation} issue, potential field methods were proposed in
\cite{yas2016excitation} and \cite{zhou2016optimal}. However, these methods
generally require global information such as target states. In this paper, we
have assumed that initial trajectories have been selected so that all
necessary events are excited. It is also possible to address this issue by
having each agent create an initial estimated potential field, until all
necessary events are excited. In addition, here we have also assumed that all
agent communications are without delays, in particular when non-local events
$\rho_{i}^{0}$ for $i\not \in \mathcal{T}_{j}(t)$ are communicated to agent
$j$ through multiple hops. The presence of delays generally affects
(\ref{eq:paramCost}). However, asynchronous versions of (\ref{eq:paramCost})
can still guarantee convergence (to the same local optima) under certain mild conditions (see \cite{bertsekas1999nonlinear}
and \cite{zhong2010asynchronous}).

\begin{algorithm}[h]
\caption{IPA-driven gradient desent for each agent}\label{IPAGradientDescent}
\begin{algorithmic}[1]
\State Initialize parameters $\bm{\theta}_j,\mathbf{w}_j$
\State Select an error tolerance $\epsilon > 0$ and a maximum number of iterations $n_0$
\State \textbf{repeat}:
\State \quad Compute the IPA gradient $\nabla_j J( \bar{\bm\theta}, \bar{\mathbf{w}})$
\State \quad Update $\bm \theta_j, \mathbf{w}_j$ using \eqref{eq:ParaUpdate}
\State \textbf{until} $\hspace{-1mm}\Vert \nabla_j J( \bar{\bm\theta}, \bar{\mathbf{w}}) \Vert \hspace{-1mm} < \hspace{-1mm}\epsilon$ or number of iterations exceeds $ n_{0}$
\State Set the optimized parameter $\bm \theta_j^{*} = \bm \theta_j, \mathbf{w}_j^{*} =  \mathbf{w}_j$.
\end{algorithmic}
\label{al:IPAGradDesc}
\end{algorithm}

\section{Simulation Examples\label{Simulation}}

We present two simulation examples to demonstrate the performance of the
decentralized scheme described in Theorem \ref{thm1}.

\begin{figure}[pb]
\centering
\includegraphics[width=1\linewidth]{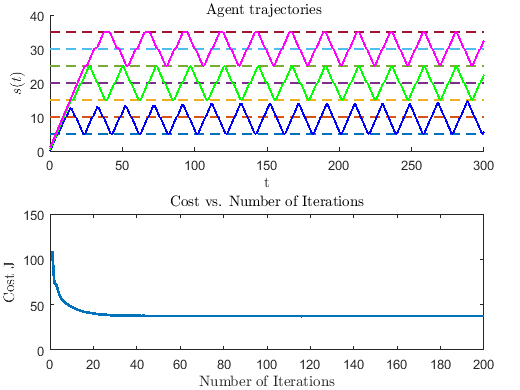}\caption{{\protect\small \textquotedblleft
Almost decentralized\textquotedblright\ optimization using Theorem \ref{thm1}.
Top plot: optimal agent trajectories. Bottom plot: cost as a function of
number of iterations with $J^{\star}=37.38$.}}%
\label{fig:Sim1_3A7T_coop}%
\end{figure}

In the first example, three homogeneous agents are allocated to persistently
monitor seven targets in the 1D mission space for $T=300$ seconds. The targets
are located at $x_{i}=5i$ for $i=1,\dots,7$. The uncertainty dynamics in
\eqref{eq:multiDynR} are defined by the parameters $A_{i}=1$, $B_{i}=5$, with
initial values $R_{i}(0)=1$ for $i=1,\ldots,7$. Each agent has a sensing range
of $r=3$ and is initialized with $s_{j}(0)=0.5(j-1)$, $u_{j}(0)=1$,
$\bm\theta_{1}^{0}=\left[  5,10,15,10,5,\ldots\right]  $, $\bm\theta_{2}%
^{0}=\left[  15,20,25,20,15,\ldots\right]  $, $\bm\theta_{3}^{0}=\left[
25,30,35,30,25,\ldots\right]  $, and $\mathbf{w}_{j}^{0}=\left[
0.5,0.5,0.5,\ldots\right]  $ for all $j=1,\dots,3$. Results of the method in
Theorem \ref{thm1} are shown in Fig. \ref{fig:Sim1_3A7T_coop}. The top plot
depicts the optimal trajectories of each agent determined after $200$
iterations of (\ref{eq:ParaUpdate}), while the bottom plot shows the overall
cost $J(\bm\theta,\mathbf{w})$ as a function of iteration number. 
The final cost is $J^{\ast}=37.38$. The
exact same results (not shown here) as in Fig. \ref{fig:Sim1_3A7T_coop} were also
obtained through the centralized scheme (\ref{CentralGradAlgo}) where all
information is available to every agent. This shows the effectiveness of the
method in Theorem \ref{thm1}.

\begin{figure}[ptb]
\centering
\includegraphics[width=0.95\linewidth]{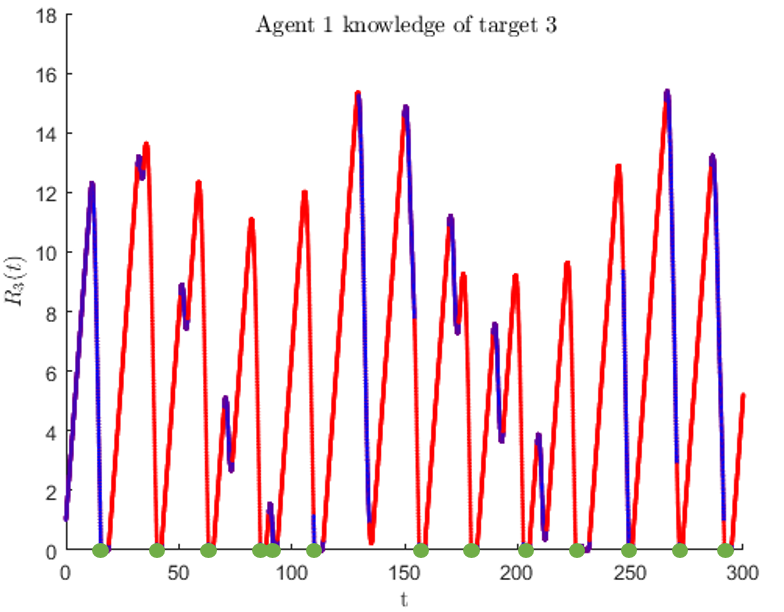}\caption{{\protect\small Red curve: $R_{3}(t)$, the state of target 3. Blue segments: $R_{3}(t)$ known to agent $1$ when its trajectory includes target 3 in its neighborhood. Green dots: instants when agent $1$ receives non-local events $\rho_{3}^{0}$.}}%
\label{fig:Sim1_RvsRest}%
\end{figure}

As pointed out earlier, the method of Theorem 1 does not involve any knowledge
by agent $j$ of the states of targets $i\not \in \mathcal{T}_{j}(t)$. This is
illustrated in Fig. \ref{fig:Sim1_RvsRest} which shows (in blue) the fraction
of time that agent 1 has any information on the state of target 3 because it
happens that $3\in$ $\mathcal{T}_{1}(t)$. The rest of the time (shown in red)
agent 1 is unable to accurately estimate the state of this target, but such
information is unnecessary. The agent only needs a small subset of is
non-local information, as illustrated by the green dots in Fig.
\ref{fig:Sim1_RvsRest}.

\begin{figure}[ptb]
\centering
\includegraphics[width=1\linewidth]{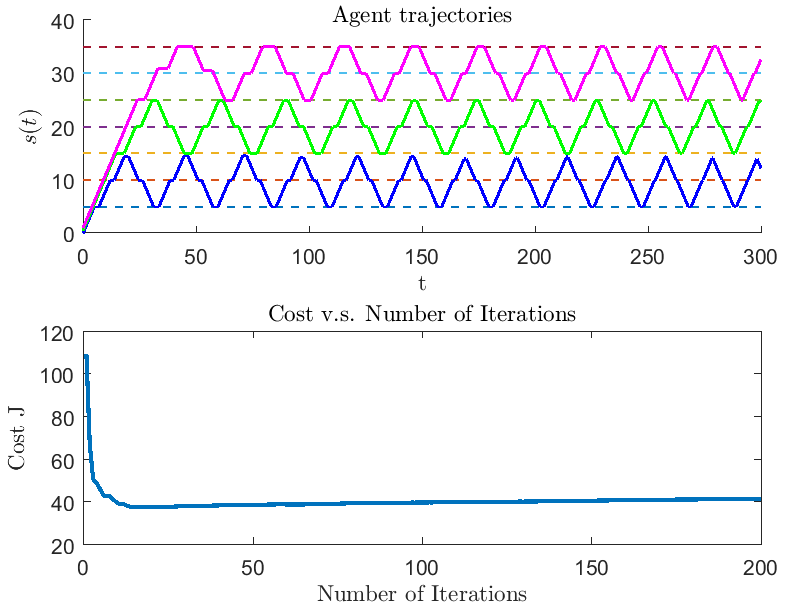}\caption{{\protect\small Fully
decentralized optimization without any non-local information. Top plot:
optimal agent trajectories. Bottom plot: cost as a function of number of
iterations with $J^{\star}=41.66$.}}%
\label{fig:Sim1_3A7T_no_global}
\end{figure}

The second example uses the same environment as the first one and agents start
with the same initial trajectories. However, we eliminate the non-local
information (condition $3$ in Theorem \ref{thm1}) and each agent calculates
its own IPA-based gradient using only local information in the set $I_{j}(t)$.
Figure \ref{fig:Sim1_3A7T_no_global} shows the results after 200 of iterations
of (\ref{eq:ParaUpdate}). Note that without non-local information, each agent
tends to spend more time dwelling on the local targets instead of better
coordinating with the other agents. Therefore, the final cost after
convergence increases from $37.38$ to $41.66$. Even though the gradient
estimate for agent $j$ is no longer accurate without the $\rho_{i}^{0}$ event
information when $i\not \in \mathcal{T}_{j}(t)$, the cost still decreases and
converges as shown in Fig. \ref{fig:Sim1_3A7T_no_global}, illustrating the
robustness of the IPA-based gradient descent method.

\section{Conclusions and Future Work}

The decentralization of multi-agent systems that involve the interaction of
agents with \textquotedblleft points of interest\textquotedblright%
\ (targets)\ in their mission space is particularly challenging. We have shown
that in 1D persistent monitoring problems an optimal centralized solution can
be recovered by an event-driven \textquotedblleft almost
decentralized\textquotedblright\ algorithm which significantly reduces
communication costs while yielding the same performance as the centralized
algorithm. In particular, each agent uses only local information except for
one event requiring communication with a non-neighbor agent when it occurs.

In addition to the event excitation issue mentioned following Theorem 1 and
incorporating communication delays in the algorithm we have developed, the
extension of this approach to the 2D case is the subject of ongoing research.
The derivations in this paper that lead to the \textquotedblleft almost
decentralized\textquotedblright\ IPA\ gradient evaluation apply to
2D-trajectories as long as these trajectories have a parametric form and the
number of parameters is finite. Moreover, the derivation holds if agents move
in straight-lines under graph-limited mobility constraints as shown in
\cite{zhou2017graphPM}. If an agent trajectory in 2D is not limited to
straight lines, the constant terms in the 1D gradient derivation in
\eqref{eq:diff_Pi} will be time-varying. The decentralization then requires
agents to share part of their trajectories with their neighboring agents when
they are in the same target neighborhood. However, this additional requirement
only involves local agent information exchanges, therefore it will not affect
the framework of our \textquotedblleft almost decentralized\textquotedblright\ solution.

\bibliographystyle{IEEEtran}
\bibliography{library}

\begin{thebibliography}{10}
\providecommand{\url}[1]{#1}
\csname url@samestyle\endcsname
\providecommand{\newblock}{\relax}
\providecommand{\bibinfo}[2]{#2}
\providecommand{\BIBentrySTDinterwordspacing}{\spaceskip=0pt\relax}
\providecommand{\BIBentryALTinterwordstretchfactor}{4}
\providecommand{\BIBentryALTinterwordspacing}{\spaceskip=\fontdimen2\font plus
\BIBentryALTinterwordstretchfactor\fontdimen3\font minus
  \fontdimen4\font\relax}
\providecommand{\BIBforeignlanguage}[2]{{%
\expandafter\ifx\csname l@#1\endcsname\relax
\typeout{** WARNING: IEEEtran.bst: No hyphenation pattern has been}%
\typeout{** loaded for the language `#1'. Using the pattern for}%
\typeout{** the default language instead.}%
\else
\language=\csname l@#1\endcsname
\fi
#2}}
\providecommand{\BIBdecl}{\relax}
\BIBdecl

\bibitem{zhong2011distributed}
M.~Zhong and C.~G. Cassandras, ``Distributed coverage control and data
  collection with mobile sensor networks,'' \emph{IEEE Trans. on Automatic
  Control}, vol.~56, no.~10, pp. 2445--2455, 2011.

\bibitem{michael2011persistent}
N.~Michael, E.~Stump, and K.~Mohta, ``Persistent surveillance with a team of
  mavs,'' in \emph{Proc. IEEE/RSJ Intl. Conf. Intelligent Robots Systems},
  2011, pp. 2708--2714.

\bibitem{leonard2010coordinated}
N.~E. Leonard, D.~A. Paley, R.~E. Davis, D.~M. Fratantoni, F.~Lekien, and
  F.~Zhang, ``Coordinated control of an underwater glider fleet in an adaptive
  ocean sampling field experiment in monterey bay,'' \emph{Journal of Field
  Robotics}, vol.~27, no.~6, pp. 718--740, 2010.

\bibitem{cassandras2013optimal}
C.~G. Cassandras, X.~Lin, and X.~Ding, ``An optimal control approach to the
  multi-agent persistent monitoring problem,'' \emph{IEEE Trans. on Automatic
  Control}, vol.~58, no.~4, pp. 947--961, 2013.

\bibitem{Smith:2012fq}
S.~L. Smith, M.~Schwager, and D.~Rus, ``{Persistent Robotic Tasks: Monitoring
  and Sweeping in Changing Environments},'' \emph{IEEE Trans. on Robotics},
  vol.~28, no.~2, pp. 410--426, Apr. 2012.

\bibitem{lin2015optimal}
X.~Lin and C.~G. Cassandras, ``An optimal control approach to the multi-agent
  persistent monitoring problem in two-dimensional spaces,'' \emph{IEEE Trans.
  on Automatic Contr.}, vol.~60, no.~6, pp. 1659--1664, 2015.

\bibitem{stump2011multi}
E.~Stump and N.~Michael, ``Multi-robot persistent surveillance planning as a
  vehicle routing problem,'' in \emph{Proc. IEEE Conf. on Automation Science
  and Engineering}, 2011, pp. 569--575.

\bibitem{yu2017optimal}
X.~Yu, S.~B. Andersson, N.~Zhou, and C.~G. Cassandras, ``Optimal dwell times
  for persistent monitoring of a finite set of targets,'' in \emph{Proc.
  American Control Conference (ACC)}, 2017, pp. 5544--5549.

\bibitem{horling2004survey}
B.~Horling and V.~Lesser, ``A survey of multi-agent organizational paradigms,''
  \emph{The Knowledge Engineering Review}, vol.~19, no.~04, pp. 281--316, 2004.

\bibitem{zhou2016optimal}
N.~Zhou, X.~Yu, S.~B. Andersson, and C.~G. Cassandras, ``Optimal event-driven
  multi-agent persistent monitoring of a finite set of targets,'' in
  \emph{Proc. IEEE Conf. on Decision and Control}, 2016, pp. 1814--1819.

\bibitem{cassandras2010perturbation}
C.~G. Cassandras, Y.~Wardi, C.~G. Panayiotou, and C.~Yao, ``Perturbation
  analysis and optimization of stochastic hybrid systems,'' \emph{European
  Journal of Control}, vol.~16, no.~6, pp. 642--661, 2010.

\bibitem{wardi2010unified}
Y.~Wardi, R.~Adams, and B.~Melamed, ``A unified approach to infinitesimal
  perturbation analysis in stochastic flow models: the single-stage case,''
  \emph{IEEE Trans. on Autom. Contr.}, vol.~55, no.~1, pp. 89--103, 2010.

\bibitem{leahy2016provably}
K.~Leahy, D.~Zhou, C.-I. Vasile, K.~Oikonomopoulos, M.~Schwager, and C.~Belta,
  ``Provably correct persistent surveillance for unmanned aerial vehicles
  subject to charging constraints,'' in \emph{Experimental Robotics}.\hskip 1em
  plus 0.5em minus 0.4em\relax Springer, 2016, pp. 605--619.

\bibitem{zhong2010asynchronous}
M.~Zhong and C.~G. Cassandras, ``Asynchronous distributed optimization with
  event-driven communication,'' \emph{IEEE Transactions on Automatic Control},
  vol.~55, no.~12, pp. 2735--2750, 2010.

\bibitem{ren2008distributed}
W.~Ren and N.~Sorensen, ``Distributed coordination architecture for multi-robot
  formation control,'' \emph{Robotics and Autonomous Systems}, vol.~56, no.~4,
  pp. 324--333, 2008.

\bibitem{olfati2007consensus}
R.~Olfati-Saber, J.~A. Fax, and R.~M. Murray, ``Consensus and cooperation in
  networked multi-agent systems,'' \emph{Proceedings of the IEEE}, vol.~95,
  no.~1, pp. 215--233, 2007.

\bibitem{bryson1975applied}
A.~E. Bryson and Y.-C. Ho, \emph{Applied optimal control: optimization,
  estimation and control}.\hskip 1em plus 0.5em minus 0.4em\relax CRC Press,
  1975.

\bibitem{yas2016excitation}
Y.~Khazaeni and C.~G. Cassandras, ``Event excitation for event-driven control
  and optimization of multi-agent systems,'' in \emph{Proc. IEEE Intl. Workshop
  on Discrete Event Systems}, 2016, pp. 197--202.

\bibitem{bertsekas1999nonlinear}
D.~P. Bertsekas, \emph{Nonlinear programming}.\hskip 1em plus 0.5em minus
  0.4em\relax Athena scientific Belmont, 1999.

\bibitem{zhou2017graphPM}
N.~Zhou, C.~G. Cassandras, X.~Yu, and S.~B. Andersson, ``Optimal event-driven
  multi-agent persistent monitoring with graph-limited mobility,'' in
  \emph{Proc. of 20th IFAC World Congress}, 2017, pp. 2217--2222.

\end{thebibliography}

\end{document}